\documentclass[1p,number,times]{elsarticle}
\usepackage[pagewise]{lineno}
\usepackage{amsfonts}
\usepackage{amssymb}
\usepackage{mathrsfs}
\usepackage{amsmath}
\usepackage{enumerate}
\usepackage{graphics}
\usepackage{mathtools}
\usepackage{microtype}
\usepackage[defaultlines=4,all]{nowidow}
\usepackage{xcolor}
\usepackage[czech, english]{babel}
\newcommand{\ad}{\overline{\rm d}}
\newcommand{\bomega}{\ensuremath{{\rm b}\omega}}
\newcommand{\ERCagreement}{{\begin{minipage}[t]{.72\textwidth}\vspace{-5pt} This paper is part of a project that has received funding from the European Research Council (ERC) under the European Union's Horizon 2020 research and innovation programme (grant agreement No 810115 -- {\sc Dynasnet}). \end{minipage}\hfill\begin{minipage}[t]{.2\textwidth}
\vspace{-6pt}
{\includegraphics[width=\textwidth]{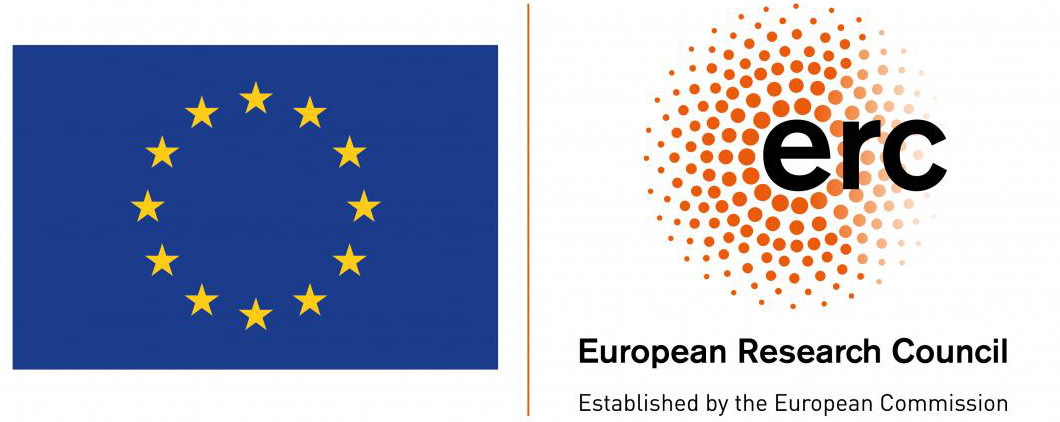}}
\end{minipage}\hfill\hfill}}

\makeatletter
\long\def\MaketitleBox{%
  \resetTitleCounters
  \def\baselinestretch{1}%
  \begin{\elsarticletitlealign}%
   \def\baselinestretch{1}%
    \Large\@title\par\vskip18pt
  \ifdoubleblind
    \vspace*{2pc}
  \else
    \normalsize\elsauthors\par\vskip10pt
    \footnotesize\itshape\elsaddress\par\vskip36pt
  \fi
  	\begin{center}
  		{\itshape \dedicated}
  	\end{center}
    \hrule\vskip12pt
    \ifvoid\absbox\else\unvbox\absbox\par\vskip10pt\fi
    \ifvoid\keybox\else\unvbox\keybox\par\vskip10pt\fi
    \hrule\vskip12pt
    \end{\elsarticletitlealign}%
}
\makeatother

\newtheorem{theorem}{Theorem}
\newtheorem{corollary}{Corollary}
\newtheorem{definition}{Definition}
\newtheorem{lemma}{Lemma}
\newtheorem{proposition}{Proposition}

\newtheorem{remark}{Remark}

\newtheorem{problem}{Problem}

\newtheorem{conjecture}[problem]{Conjecture}

\newenvironment{absolutelynopagebreak}
  {\par\nobreak\vfil\penalty0\vfilneg
   \vtop\bgroup}
  {\par\xdef\tpd{\the\prevdepth}\egroup
   \prevdepth=\tpd}

\newproof{proof}{Proof}
\newproof{potweak}{Proof of Theorem~\ref{thm:weak}}
\newproof{proofk2}{Proof of Theorem~\ref{thm:k2}}
\newcommand{\myqed}{\qed}
\journal{some journal}
\begin{document}
\begin{frontmatter}

\title{From $\chi$- to $\chi_p$-bounded classes\tnoteref{ERC}}
\tnotetext[ERC]{\ERCagreement}
\author{Yiting Jiang}
\ead{yjiang@irif.fr}
\address{Universit\'e de Paris, CNRS, IRIF, F-75006, Paris, France\\and Department of Mathematics, Zhejiang Normal University, China}
\author{Jaroslav Ne\v set\v ril\corref{cor1}}
\ead{nesetril@iuuk.mff.cuni.cz}
\address{Computer Science Institute of Charles University (IUUK), Praha, Czech Republic}
\author{Patrice Ossona de Mendez} 
\ead{pom@ehess.fr}
\address{Centre d'Analyse et de Math\'ematiques Sociales (CNRS, UMR 8557),
Paris, France\\
  and     Computer Science Institute of Charles University,
  Praha, Czech Republic}
\cortext[cor1]{Corresponding author}
\def\dedicated{To the memory of Robin Thomas}
\begin{abstract}
$\chi$-bounded classes are studied here in the context of star colorings and, more generally, $\chi_p$-colorings.
This fits to a general scheme of sparsity and leads to natural extensions of the notion of bounded expansion class. In this paper we solve two conjectures related to star coloring ({i.e.} $\chi_2$) boundedness. One of the conjectures is disproved and in fact we determine which weakening holds true. $\chi_p$-boundedness leads to more stability and we give structural characterizations of (strong and weak) $\chi_p$-bounded classes.
We also generalize a result of Wood relating the chromatic number of a graph to the star chromatic number of its $1$-subdivision.
As an application of our characterizations, among other things, we show that for every odd integer $g>3$ even hole-free graphs $G$ contain at most $\varphi(g,\omega(G))\,|G|$ holes of length $g$.
\end{abstract}
\end{frontmatter}

\section{Introduction}
The concept of $\chi$-boundedness was introduced by Gy\'arf\'as in 1985 in his seminal paper \cite{gyarfas1985problems}. A family of graphs $\mathscr C$ is {\em $\chi$-bound} (or {\em $\chi$-bounded}) with {\em binding function} $f$ if $\chi(H)\leq f(\omega(H))$ holds whenever $G\in\mathscr C$ and $H$ is an induced subgraph of $G$. The notion of $\chi$-boundedness has attracted much attention and motivated important conjectures (see survey \cite{scott2018survey}). Because the definition of $\chi$-boundedness involves all the induced subgraphs of the graphs in the family, it will be natural to restrict our attention to {\em hereditary} classes of graphs, that is to classes of graphs closed under induced subgraphs.

In this setting, probably the most important open conjecture is  the next one.

\begin{conjecture}[Gy\'arf\'as \cite{Gyarfas1985}, Sumner \cite{sumner1981subtrees}]
\label{conj:GyarfasSumner}
	For every tree $T$, the class of all graphs excluding $T$ as an induced subgraph is $\chi$-bounded.
\end{conjecture}

Indeed, as there exist graphs with arbitrary high girth and chromatic number \cite{ErdH1959}, excluding an induced subgraph with a cycle does not allow to bind the chromatic number by a function of the clique number. A natural alternative is to forbid some fixed graph as an induced subdivision (that is to forbid all the subdivisions of some fixed graph as induced subgraphs). This motivated the following conjecture.

\begin{conjecture}[Scott \cite{Scott1997}]
\label{conj:Scott}
	For every graph $F$ the class of all graphs excluding induced subdivisions of $F$ is $\chi$-bounded.
\end{conjecture}
This conjecture was disproved by Pawlik, Kozik, Krawczyk, Laso\'n, Micek, Trotter and Walczak \cite{Pawlik2014}.  Nevertheless the conjecture motivated several positive results and here we add to this list several new instances.
Note that for biclique-free classes of graphs (i.e. classes of graphs excluding some fixed biclique $K_{r,r}$ as a subgraph) of these conjectures hold. For Conjecture~\ref{conj:GyarfasSumner} this has been proved by Kierstead and R{\"o}dl \cite{kierstead1996applications}, while for Conjecture~\ref{conj:Scott} this follows from K{\" u}hn and Osthus \cite{kuhn2004induced}.

	Similar to the notion of $\chi$-boundedness, Karthick \cite{karthick2018star} introduced the notion of $\chi_s$-bounded class, where $\chi_s$ denotes the star chromatic number. Recall that the {\em star chromatic number} of a graph $G$, a notion introduced by Gr\" unbaum \cite{Gruenbaum1973},  is the minimum number of colors in a proper coloring of $G$ with the property that any two color classes induce a star forest. 
	In this setting,  two conjectures were proposed. As customary, by an {\em $F$-free} graph we mean a graph no containing $F$ as an induced subgraph and, for a family $\mathscr F$ of graphs, an {\em $\mathscr F$-free} graph is a graph which is $F$-free for all $F\in\mathscr F$.

\begin{conjecture}[Karthick \cite{karthick2018star}]
\label{conj:k1}
	The class of all $K_{1,t}$-free graphs (where $t\geq 3$) is $\chi_s$-bounded.
\end{conjecture}
\begin{conjecture}[Karthick \cite{karthick2018star}]
\label{conj:k2}
	For any tree $T$, the class of all $(T,C_4)$-free graphs is $\chi_s$-bounded.
\end{conjecture}

In Section~\ref{sec:chis} we prove Conjecture~\ref{conj:k1} (see Theorem~\ref{thm:k1}), disprove Conjecture~\ref{conj:k2}, and characterize those classes of $(T,K_{r,t})$-free graphs (with $T$ a forest) that are $\chi_s$-bounded.

\begin{theorem}
\label{thm:k2}
		Let $T$ be a forest and let $r\leq t$ be positive integers. Then the class $\mathscr C$ of 
	all $(T,K_{r,t})$-free graphs is $\chi_s$-bounded if and only if $r=1$ or $T$ is a subgraph of the $1$-subdivision of a tree. 
\end{theorem}

In \cite{Taxi_tdepth}, a generalization of the chromatic number was proposed, which defines a non-decreasing sequence $\chi_1,\chi_2,\dots$ of graph invariants, where $\chi_1$ is the usual chromatic number (i.e. $\chi_1=\chi$), $\chi_2$ is the star chromatic number (i.e. $\chi_2=\chi_s$), and $\chi_p$ is the minimum number of colors of a low tree-depth coloring with parameter $p$ (see Section~\ref{sec:prelim}).

The notion of bounded expansion captures uniform sparsity of graph classes. Formally, a class $\mathscr C$ has {\em bounded expansion} if the shallow minors at depth $r$ of graphs in $\mathscr C$ have their average degree bounded by some function of $r$ (see Section~\ref{sec:prelim}). 
Some characterizations of bounded expansion classes will be of prime importance here and we review them now. It is one of the important features of the theory of sparsity that classes with bounded expansion can be characterized in many different ways.

For a graph $G$ and a non-negative integer $r$ we denote by ${\rm TM}_r(G)$ the class of all graphs $H$, with the property that a $(\leq r)$-subdivision of $H$ (i.e. a graph obtained from $H$ by subdividing each edge by at most $r$ vertices) is a subgraph of $G$. (Such a graph is also a \emph{shallow topological minor} of $G$ at depth $r/2$.) More generally, if $\mathscr C$ is a class of graphs we define ${\rm TM}_r(\mathscr C)=\bigcup_{G\in\mathscr C}{\rm TM}_r(G)$. 
Also, following \cite{DVORAK2018143}, we denote by 
${\rm ITM}_r^e(G)$  the class of all graphs $H$ whose (exact) $r$-subdivision $H^{(r)}$ is an induced subgraph of $G$, and let 
${\rm ITM}_r^e(\mathscr C)=\bigcup_{G\in\mathscr C}{\rm ITM}_r^e(G)$.

For a graph invariant $f$ and a class $\mathscr C$ it will be convenient to define $f(\mathscr C)=\sup_{G\in\mathscr C}f(G)$.
For example, $\chi(\mathscr C)$ is the supremum of the chromatic numbers of the graphs in $\mathscr C$ and $\ad(\mathscr C)$ is  the supremum of the average degrees $\ad(G)$ of the graphs $G$ in $\mathscr C$.
Classes with bounded expansion are characterized by means of the average degrees of topological minors.
\begin{lemma}[{\cite[Theorem 11]{Dvov2007}}]
\label{lem:BEtm}
	A class $\mathscr C$ has bounded expansion if and only if for every non-negative integer $r$ we have
	\[
		\ad({\rm TM}_r(\mathscr C))<\infty,
	\]
	where $\ad$ denotes the average degree.
\end{lemma}
Bounded expansion classes are also characterized by means of the $\chi_p$-invariants.
\begin{lemma}[{\cite[Theorem 7.1]{POMNI}}]
\label{lem:BEchip}
		A class $\mathscr C$ has bounded expansion if and only if for every positive integer $p$ we have
	\[
		\chi_p(\mathscr C)<\infty.
	\]
\end{lemma}
(See \cite{Sparsity} for a general background of sparsity.)

Motivated by these characterizations, we consider in this paper two generalizations of the notion of $\chi$-boundedness.
\begin{definition}
	A hereditary class $\mathscr C$ is {\em strongly $\chi_p$-bounded} if, 
	for every $G\in\mathscr C$ we have $\chi_p\leq f_p(\omega(G))$ (for some fixed binding function $f_p$)
\end{definition}

\begin{definition}
	A hereditary class $\mathscr C$ is {\em weakly $\chi_p$-bounded} if, for every $G\in\mathscr C$ we have $\chi_p\leq g_p(\omega({\rm TM}_{p-1}(G)))$ (for some fixed binding function $g_p$).
\end{definition}

This second definition may look arbitrary at first glance. 
The reason why we consider ${\rm TM}_{p-1}(G)$ in the definition of weakly $\chi_p$-bounded classes is that $\chi_p(G)$ has a lower bound in terms of $\omega({\rm TM}_{p-1}(G))$ (see Lemma~\ref{lem:chiptm})  but not in terms of $\omega({\rm TM}_{p}(G))$. To see this, let $G^{(p)}$ denote the \emph{$p$-subdivision} of a graph $G$, that is the graph obtained by replacing each edge of $G$ by a path of length $p+1$.
 Then we have $\chi_p(K_n^{(p)})=p+1$ and $\omega({\rm TM}_{p}(K_n^{(p)}))=n$, while $\omega({\rm TM}_{p-1}(K_n^{(p)}))=2$.
This lower bounds suggest a  generalization 
of some inequalities \cite{wood2005acyclic} binding the chromatic number $\chi(G)$ of a graph $G$ and the star chromatic number $\chi_2(G^{(1)})$ of the $1$-subdivision of $G$  into inequalities binding $\chi(G)$ and the $p$-th chromatic number $\chi_p(G^{(p-1)})$ of the $(p-1)$-subdivision of $G$ (see Theorem~\ref{thm:chip-sub} in Section~\ref{sec:weak}).

In Section~\ref{sec:strong}
we give several examples of strongly $\chi_p$-bounded classes, including induced subgraphs of the $d$-power of graphs in a bounded expansion class, claw-free graphs, trivially perfect graphs, even hole-free graphs, and split graphs, and then 
 we give a characterization of strongly $\chi_p$-bounded classes (where undefined notions will be defined in Section~\ref{sec:strong}).

Strongly $\chi_p$-bounded  classes are structurally characterized by the following result.
\begin{theorem}
\label{thm:strong}
	Let $\mathscr C$ be a hereditary class of graphs.
	Then the following are equivalent:
\begin{enumerate}[\rm (i)]
	\item\label{enum:m1} The class $\mathscr C$ is strongly $\chi_p$-bounded for every integer $p$;
	\item\label{enum:m2} For each positive integer $t$, the class $\mathscr C_t=\{G\in\mathscr C\mid\omega(G)\leq t\}$ has bounded expansion;
	\item\label{enum:m3} 
	The class $\mathscr C$ has  $\omega$-bounded expansion, meaning that
	for every non-negative integer $r$ there is a function $f_r$ such that for every $G\in\mathscr C$ we have $\ad({\rm TM}_r(G))\leq f_r(\omega(G))$ (see Section~\ref{sec:strong});
	\item\label{enum:m4} The class $\mathscr C$ is $\chi$-bounded, 
	does not contain all complete bipartite graphs, and for every positive integer $r$ we have $\ad({\rm ITM}_r^e(\mathscr C))<\infty$;
	\item\label{enum:m5} Every connected acyclically oriented graph has a restricted dual for the
	 class of all orientations of graphs in $\mathscr C$.
\end{enumerate}
\end{theorem}

In Section~\ref{sec:weak}
we show that the class of complete bipartite graphs is weakly $\chi_p$-bounded but not strongly $\chi_p$-bounded.
We then give the next structural characterization of weakly $\chi_p$-bounded classes and deduce (Proposition~\ref{prop:SBEchip}) 
that first-order transductions of bounded expansion classes are weakly $\chi_p$-bounded for every $p$.
\begin{theorem}
\label{thm:weak}
	Let $\mathscr C$ be a hereditary class of graphs.
	Then the following are equivalent:
\begin{enumerate}[\rm (i)] 
		\item\label{enum:w1} The class $\mathscr C$ is weakly $\chi_p$-bounded for every positive integer $p$;
		\item\label{enum:w2} the class $\mathscr C$ and all the classes ${\rm ITM}_r^e(\mathscr C)$ ($r\geq 1$) are $\chi$-bounded;
		\item\label{enum:w3} the class $\mathscr C$ is $\chi$-bounded and for every positive integer $r$ we have $\ad({\rm ITM}_r^e(\mathscr C))<\infty$;
		\item\label{enum:w4} $\mathscr C$ is $\chi$-bounded and for each positive integer $p$ there is a function $f_p$ such that for every graph $G\in\mathscr C$ we have
			$\chi_p(G)\leq f_p(\bomega(G))$, where $\bomega(G)=\max\{s\mid K_{s,s}\subseteq G\}$;
		\item\label{enum:w5} $\mathscr C$ is $\chi$-bounded and for each positive integer $s$ the class $\{G\in\mathscr C\mid K_{s,s}\not\subseteq_i G\}$ is strongly $\chi_p$-bounded for every positive integer $p$. 
\end{enumerate}
\end{theorem}

In Section~\ref{sec:weak} we give examples of weakly $\chi_p$-bounded classes of graphs, including classes with low twin-width covers and proper vertex-minor-closed classes.

In Section~\ref{sec:appli} we give some applications. Among other things, we prove that for every odd integer $g>3$ even hole-free graphs $G$ contain at most $\varphi(g,\omega(G))\,|G|$ holes of length $g$ (Theorem~\ref{thm:holes}).

\section{Definitions and notations}
\label{sec:prelim}
We denote by $H\subseteq G$ the property that $H$ is a subgraph of $G$, and by $H\subseteq_i G$ the property that $H$ is an induced subgraph of $G$. A class of graphs is {\em monotone} if it is closed under subgraphs; it is \emph{hereditary} if it is closed under induced subgraphs.

We denote by $\omega(G)$ the {\em clique number} of $G$, i.e. $\max\{t\mid K_t\subseteq G\}$ and we define the bipartite analog, the \emph{biclique number} $\bomega(G)$ of $G$:
\[
	\bomega(G)=\max\{r\mid K_{r,r}\subseteq G\}.
\]

Note that obviously $\bomega(G)\geq\lfloor\omega(G)/2\rfloor$.
A class $\mathscr C$ with $\bomega(\mathscr C)<\infty$ is said to be {\em weakly sparse} \cite{Jiang2020} or {\em biclique free}.

A graph $H$ is a {\em $(\leq r)$-subdivision} (resp. the {\em $r$-subdivision}) of a graph $G$ if it can be obtained from $G$ by subdividing each edge by at most $r$ vertices (resp. by exactly $r$ vertices).
The $r$-subdivision of a graph $G$ is denoted by $G^{(r)}$.

Recall that 
the class ${\rm TM}_r(G)$ is the class of all graphs $H$, such that some $(\leq r)$-subdivision of $H$ is a subgraph of $G$, and that
the class ${\rm ITM}_r^e(G)$ is the class of all graphs $H$, such that the $r$-subdivision of $H$ is an induced  subgraph of $G$. 

The {\em tree-depth} of a graph $G$, denoted by ${\rm td}(G)$, is the minimum height of a rooted forest $Y$ such that $G$ is a subgraph of the closure of $Y$. Equivalently, the tree-depth of $G$ is the minimum clique number of a trivially perfect supergraph of $G$. 

The {\em tree-depth chromatic number of rank $p$} of $G$, denoted by $\chi_p(G)$ is the minimum number of colors in a vertex coloring of $G$ such that every subset $I$ of at most $p$ colors induce a subgraph $G_I$ of $G$ with tree-depth at most $|I|$. In particular, for every graph $G$ we have $\chi_1(G)=\chi(G), \chi_2(G)=\chi_s(G)$, and 
$\chi_1(G)\leq \chi_2(G)\leq\dots\leq \chi_{|G|}(G)={\rm td}(G)$.
For basic properties of tree-depth and $\chi_p$ we refer the reader to \cite{Sparsity}.

\section{$\chi_s$-bounded classes of graphs}
\label{sec:chis}

 Conjecture~\ref{conj:k1} can be proved in an easy way. 
\begin{theorem}
\label{thm:k1}
	The class $K_{1,t}$-free graphs is polynomially $\chi_s$-bounded. Precisely, the $K_{1,t}$-free graphs $G$ satisfy
\[
	\chi_s(G)=O(\omega(G)^{\frac{3(t-1)}{2}}).
\] 
\end{theorem}
\begin{proof}
	Excluding $K_{1,t}$ we get a class where the maximum degree $\Delta(G)$ of a graph $G$ is bounded by a function of its  clique number $\omega(G)$.
	Precisely we have:
	\[
		\Delta(G)< R(\omega(G),t)\leq 
		\binom{\omega(G)+t-2}{t-1}
		=O(\omega(G)^{t-1}),
	\]
	 where $R(n,m)$ is the Ramsey number for two coloring of edges of a complete graph.
	Indeed, if the neighborhood of any vertex $v$ has size $R(\omega(G),t)$ and does not contain an independent set  of size $t$ then it contains a clique  of size $\omega(G)$, a contradiction.
	Then the theorem follows from the fact that
 graphs with degree at most $D$ have star chromatic number $O(D^{3/2})$ \cite{Fertin2001}.
\myqed\end{proof} 
Let us remark that, 
by the same argument, we get that $\ad({\rm TM}_r(G))$ is bounded by some fixed function of $\omega(G)$ and $r$ for every $K_{1,t}$-free graph $G$.
In other words, the class of  $K_{1,t}$-free graphs has \emph{$\omega$-bounded expansion} (see formal definition in Section~\ref{sec:strong}). This stronger property actually holds for the more general class of all even hole-free graphs that exclude induced subdivisions of a fixed complete bipartite graph $K_{s,t}$.

We now disprove Conjecture~\ref{conj:k2} in the following more general form.
To this end, we make use of the following result of Wood.
\begin{lemma}[{\cite[Theorem 2]{wood2005acyclic}}]
\label{lem:chis-sub}
For every graph $G$ the star chromatic number of $G^{(1)}$ satisfies:
\[\sqrt{\chi(G)}\leq \chi_s(G^{(1)})\leq\max(\chi(G),3)\]
\end{lemma}


\begin{lemma}
\label{lem:notk2}
	Let $T$ be a forest that is not a subgraph of a $1$-subdivided tree, and let $k$ be a positive integer.
	Then the class of all $(T,C_4,\dots,C_{2k}$, $\text{odd hole})$-free graphs is not $\chi_s$-bounded.
\end{lemma}
\begin{proof}
	 Let $\mathscr C$ be the hereditary closure of the class of $1$-subdivisions of all the graphs with girth at least $k+1$.
	  The graphs in $\mathscr C$ are $(T,C_4,\dots,C_{2k}$, $\text{odd hole})$-free and have no triangle. 
	  It is well known that graphs of girth at least $k+1$ have unbounded chromatic number \cite{ErdH1959}. According to Lemma~\ref{lem:chis-sub} we have $\chi_s(G^{(1)})\geq\sqrt{\chi(G)}$ and hence the class $\mathscr C$ is not $\chi_s$-bounded.
\myqed\end{proof}

We now characterize those forests $T$ and those complete bipartite graphs $K_{r,t}$ with the property that the class of $(T,K_{r,t})$-free graphs is $\chi_s$-bounded. 

For this we shall need the following results, which we restate using our definitions and notations.
\begin{lemma}[{\cite[Theorem 2]{kuhn2004induced}}]
\label{lem:dense1sub}
	For all $k,r\in\mathbb N$ there exists $d=d(r,k)$ such that
	for every  graph $G$ with $K_{r,r}\not\subseteq G$ and $\ad(G)\geq d$ we have 
	\[\ad({\rm ITM}_1^e(G))\geq k\] (i.e. $G$ 
	contains an induced subdivision of a graph $H$, where the average degree of $H$ is at least $k$ and every edge of $H$ is subdivided exactly once).
\end{lemma}

\begin{lemma}[{\cite[Lemma 9]{DVORAK2018143}}]
\label{lem:zd}
Let $r,k,t\geq 1$ be integers. Let $r_0=\max(r,2^{25},t+1,tk/2)$ and $d_{r,k,t}=\frac{r_0^{11}(tk+1)}{2^6}$. Then for every graph $G$ with 
	$\max_{H\subseteq G}\ad(H)\leq t$ and 
	$\ad({\rm ITM}_k^e(G))<r$	
	we have
	\[\ad({\rm TM}_k(G))<\ad({\rm TM}_{k-1}(G))+d_{r,k,t}.\]
\end{lemma} 

We shall need the following nice results of Dvo\v r\'ak.
\begin{lemma}[{\cite[Lemma 7]{Dvov2007}}]
\label{lem:1subchr}
For every integer $c$ there exists an integer $d=d(c)$ such that every graph with average degree at least $d$ contains as a subgraph the $1$-subdivision of a graph with chromatic number $c$.
\end{lemma}
\begin{lemma}[{\cite[Corollary 4]{Dvov2007}}]
\label{lem:chichis}
	Let $\mathscr C$ be any class of graphs with $\chi(\mathscr C)<\infty$. Then
	$\chi_s(\mathscr C)<\infty$ if and only if there exists $c$ 
 such that if $G^{(1)}$ 
 is a subgraph of a graph in $\mathscr C$
then $\chi(G)\leq c$.
\end{lemma}
Lemma~\ref{lem:chichis} can be equivalently restated by the fact that $\chi_s(G)$ and $\chi({\rm TM}_1(G))$ are bound to each other, in the sense that there exist functions $f$ and $g$ such that $\chi_s(G)\leq f(\chi({\rm TM}_1(G)))$ and $\chi({\rm TM}_1(G))\leq g(\chi_s(G))$.

We are now in the position to prove Theorem~\ref{thm:k2} stated in the introduction.

\begin{proofk2}
	If $r=1$ then the result follows from Theorem~\ref{thm:k1}. Thus we can assume $r\geq 2$.
	
	According to Lemma~\ref{lem:notk2}, if $\mathscr C$ is $\chi_s$-bounded then $T$ is a subgraph of a $1$-subdivided tree. 
	Assume $T$ is a subgraph of a $1$-subdivided tree $T'$. We can assume that $T$ is an induced subgraph of $T'$. We consider the class $\mathscr C'$ of all $(T',K_{r,r})$-free graphs (which includes $\mathscr C$). Let $k$ be an integer and let $\mathscr C_k'$ be the subclass of $\mathscr C'$ of graphs with clique number at most $k$.  By an easy Ramsey argument, the graphs in the class $\mathscr C'$ exclude some $K_{r',r'}$ as a (non induced) subgraph. Note that $\mathscr C_k'$ is a hereditary class.
	
	Assume the average degree of graphs in $\mathscr C_k'$ is arbitrarily large. By Lemma~\ref{lem:dense1sub} (and as $\mathscr C_k'$ is hereditary), the class 
	 $\mathscr C_k'$ contains  $1$-subdivided graphs with arbitrarily large average degree, in which we can easily find a copy of $T'$, contradicting our hypothesis. 
	
	Thus the average degree of the graphs in $\mathscr C_k'$ is bounded by some constant $C(k)$, and so is the chromatic number. Assume that the graphs in $\mathscr C_k'$ have unbounded $\chi_s$. Then, according to Lemma~\ref{lem:chichis} we can find graphs $G$ in $\mathscr C_k'$ such that the average degree of the graphs in ${\rm TM}_1(G)$ is arbitrarily large. Then, according to Lemma~\ref{lem:zd} the average degree of the graphs in ${\rm ITM}_1^e(G)$ are also arbitrarily large, again leading to a contradiction.
\myqed\end{proofk2}

\section{Strongly $\chi_p$-bounded classes of graphs}
\label{sec:strong}
Theorem~\ref{thm:strong} characterizes classes that are strongly $\chi_p$-bounded. It
will directly follow from the following results stated below as Lemma~\ref{lem:OBE0} (equivalence of \eqref{enum:m1},\eqref{enum:m2}, and \eqref{enum:m3}), Lemma~\ref{lem:OBEi} (equivalence of \eqref{enum:m3} and \eqref{enum:m4}), and Lemma~\ref{lem:OBErd} (equivalence of \eqref{enum:m3} and \eqref{enum:m5}).

	A class $\mathscr C$ has {\em $\omega$-bounded expansion} if there exists a function $f:\mathbb N\times\mathbb N\rightarrow\mathbb R$ such that for every $G\in\mathscr C$ and every non-negative integer $r$ we have
	\begin{equation}
		\ad({\rm TM}_r(G))\leq f(\omega(G),r).
	\end{equation}

\begin{lemma}
\label{lem:OBE0}
	For a class $\mathscr C$ the following are equivalent:
	\begin{enumerate}[\rm (i)]
		\item\label{enum:1} for each integer $t$ the subclass $\mathscr C_t$ of all the $K_t$-free graphs in $\mathscr C$ has bounded expansion;
		\item\label{enum:2} the class $\mathscr C$ has $\omega$-bounded expansion;
		\item\label{enum:3} the class $\mathscr C$ is strongly $\chi_p$-bounded for each integer $p$. Explicitly, for every integer $p$ there exists a function $f_p$ such that $\chi_p(G)\leq f_p(\omega(G))$ for every graph $G$ in the class $\mathscr C$.
	\end{enumerate} 
\end{lemma}
\begin{proof}
	$\eqref{enum:1}\Leftrightarrow\eqref{enum:2}$:
	Assume \eqref{enum:1}. Then, for each integer $t$, there exists a function $g_t:\mathbb N\rightarrow\mathbb N$ such that 
		for every shallow minor $H$ of $G$ at depth $r$ we have $\|H\|/|H|\leq g_t(r)$. Defining $f(t,r)=g_t(r)$ we deduce that $\mathscr C$ has $\omega$-bounded expansion. 
	The converse implication is also obvious.
			
	$\eqref{enum:2}\Leftrightarrow\eqref{enum:3}$:
	Assume~\eqref{enum:2}.
According to Lemma~\ref{lem:BEchip}, for every integer $p$ there is a constant $c_t(p)$ with $\chi_p(G)\leq c_t(p)$ for every $G\in\mathscr C_t$. Hence, defining, $f_p(t)=c_t(p)$ we get $\chi_p(G)\leq f_p(\omega(G))$ for every $G\in\mathscr C$.
	Conversely, assume~\eqref{enum:3}.
 Then $\chi_p(G)$ is bounded by the constant $f_p(t)$ on $\mathscr C_t$. Thus, according to Lemma~\ref{lem:BEchip}, $\mathscr C_t$ has bounded expansion. Hence $\mathscr C$ has $\omega$-bounded expansion.
\myqed\end{proof}

It should be noticed that function $f_p$ appearing in Item~(iii) of Lemma~\ref{lem:OBE0} can be bounded in terms of the function $f_1$ and the diagonal terms $f_p(p)$. This is  a direct corollary of the next proposition.

\begin{proposition}
Let $\mathscr C$ be a hereditary strongly $\chi_p$-bounded class and let
$a_p=\max\{\chi_p(G): G\in\mathscr C\text{ and }\omega(G)\leq p\}$.
Then, for every graph $G\in\mathscr C$ and every positive integer $p$ we have
\[
	\chi_p(G)\leq \chi(G)\,a_p^{\binom{\chi(G)-1}{p-1}}.
\]
\end{proposition}
\begin{proof}
	Let $p$ be a positive integer, let $G\in\mathscr C$, let $\chi=\chi(G)$, and let $c:V(G)\rightarrow [\chi]$ be a proper coloring of $G$.
	For any subset $I$ of $p$ colors in $[\chi]$, let $G_I$ be the subgraph of $G$ induced by vertices with color in $I$.
	Note that $\omega(G_I)\leq p$ and $G_I\in\mathscr C$ as $\mathscr C$ is hereditary. Thus there exists a $\chi_p$-coloring of $G_I$ with $a_p$ colors. Let $\gamma_I:V(G_I)\rightarrow[a_p]$ be such a $\chi_p$-coloring.
	For $v\in V(G)$ define $g_v:\binom{[\chi]\setminus\{c(v)\}}{p-1}\rightarrow [a_p]$ by $g_v(J)=\gamma_{J\cup\{c(v)\}}(v)$. Consider the coloring $\zeta: v\mapsto \zeta(v)=(c(v),g_v))$. This coloring uses at most
	$\chi\, a_p^{\binom{\chi-1}{p-1}}$ colors. We now prove that $\zeta$ is a $\chi_p$-coloring of $G$. Let $\zeta_1=(c_1,g_1),\dots,\zeta_p=(c_p,g_p)$ be $p$ $\zeta$-values,	and let $I$ be a subset of size $p$ of $[\chi]$ that includes $c_1,\dots,c_p$. For $1\leq i\leq p$ let $V_i$ be the set of vertices $v$ of $G$ with $\zeta(v)=\zeta_i$. Obviously, $\bigcup_{i\in I}V_i\subseteq V(G_I)$. Let $1\leq i\leq p$ and let $v\in V_i$. Then $\gamma_I(v)=g_v(I\setminus\{c(v)\})=g_i(I\setminus\{c_i\})$. Thus all the vertices in $V_i$ have the same $\gamma_I$-color. It follows that the subgraph of $G$ induced by $\zeta$-colors $\zeta_1,\dots,\zeta_p$ is an induced subgraph of $G_I$ induced by at most $p$ $\gamma_I$-colors thus has tree-depth at most $p$. Hence $\zeta$ is a $\chi_p$-coloring of $G$.
\myqed\end{proof}

Inspired by (non valid) Conjecture~\ref{conj:Scott}, let us mention the following "positive" result.
\begin{lemma}[{\cite[Theorem 4]{DVORAK2018143}}]
\label{lem:dvo-sub}
	For every graph $H$ and a positive integer $r$, if $\mathscr C$ is a class of graphs that do not contain $K_r$, $K_{r,r}$, and any subdivision of $H$ as an induced subgraph, then $\mathscr C$ has bounded expansion.
\end{lemma}
This has the following immediate consequence.
\begin{corollary}
	For every graph $H$ and a positive integer $r$, the class of all graphs excluding both an induced $K_{r,r}$ and all induced subdivisions of $H$ has $\omega$-bounded expansion.
\end{corollary}

Consequently the following classes have $\omega$-bounded expansion (hence are $\chi$-bounded, $\chi_s$-bounded and, more generally, $\chi_p$-bounded):
\begin{itemize}
	\item Any class of graphs with bounded stability number. Indeed, if $\alpha(G)<t$ then $G$ excludes induced $K_{t,t}$ and (every induced subdivision of) $tK_1$.
	This includes, for instance, the class of the complements of shift-graphs.
	\item Any class of graphs such that the neighborhood of every vertex has bounded stability number. (Indeed, this boils down to excluding $K_{1,t}$.) Note that this includes claw-free graphs.
	\item Any hereditary class of graphs excluding a complete bipartite graph and having a bound on the diameter. This includes trivially perfect graphs.
	\item The class of all even hole-free graphs. Actually, it is sufficient to consider theta-free graphs, as this amounts to exclude all subdivisions of $K_{2,3}$. More generally, graphs excluding all induced subdivisions of $K_{r,t}$ form an $\omega$-bounded expansion class.
	\item The class of all split graphs. Indeed, split-graphs exclude $C_4,C_5$, and $\overline{C}_4=2K_2$. Note that every subdivision of $2K_2$ includes $2K_2$ as an induced subgraph.
\end{itemize}

Some further examples of classes with $\omega$-bounded expansion (and thus strongly $\chi_p$-bounded) can be obtained from classes with bounded expansion. It follows from \cite{power} that if a class $\mathscr C$ has bounded expansion and $d$ is a positive integer, then the class $\mathscr C^d=\{G^d\mid G\in\mathscr C\}$ has $\omega$-bounded expansion. Here $G^d$ denotes the {\em $d$-th power} of $G$, that is the graph with vertex set $V(G)$ in which two vertices are adjacent if their distance in $G$ is at most $d$. These examples are typical and they lead to the following:

\begin{lemma}
\label{lem:OBEi}
	Let $\mathscr C$ be a hereditary class of graphs. Then $\mathscr C$ has $\omega$-bounded expansion if and only if
	\begin{enumerate}[\rm (i)]
		\item \label{enum:OBE1} $\mathscr C$ is $\chi$-bounded,
		\item \label{enum:OBE2} $\mathscr C$ does not contain all complete bipartite graphs, and
		\item \label{enum:OBE3} for every integer $r\geq 1$
		we have $\ad({\rm ITM}_r^e(\mathscr C))<\infty$.
	\end{enumerate}
\end{lemma}
\begin{proof}
	Assume the conditions \eqref{enum:OBE1},\eqref{enum:OBE2}, and \eqref{enum:OBE3} are satisfied. 
	Let $t\in\mathbb N$ and let $\mathscr C_t=\{G\in\mathscr C\mid \omega(G)\leq t\}$. By \eqref{enum:OBE2} graphs in $\mathscr C_t$ exclude some $K_{s,s}$ as an induced subgraph. By Ramsey theorem, as they also exclude $K_t$ they exclude $K_{R(t,s),R(t,s)}$ as a subgraph. According to \eqref{enum:OBE3}, there exists a constant $d=\ad({\rm ITM}_1^e(\mathscr C))$ such that if the $1$-subdivision of a graph $H$ is an induced subgraph of a graph $G\in\mathscr C_t$ then the average degree of $H$ is at most $d$.
	From this property and the exclusion of $K_{R(t,s),R(t,s)}$ as a subgraph we deduce, according to Lemma~\ref{lem:dense1sub} that there exists a constant $d'$ such that every graph in $\mathscr C_t$ (as well as every induced subgraph of graphs in $\mathscr C_t$ as $\mathscr C_t$ is hereditary) has average degree at most $d'$. Then, according to Lemma~\ref{lem:zd}, it follows from this property, \eqref{enum:OBE1} and \eqref{enum:OBE3} that $\mathscr C$ has $\omega$-bounded expansion.
	
	Assume that $\mathscr C$ has $\omega$-bounded expansion. Conditions \eqref{enum:OBE1} and \eqref{enum:OBE2} are obviously satisfied. Let $\mathscr C_3$ be the subclass of all triangle-free graphs in $\mathscr C$. The class $\mathscr C_3$ has bounded expansion hence for every integer $r$ there is a constant $f(r)$ such that if the $r$-subdivision $G$ of a graph $H$ belongs to $\mathscr C_3$ then $\ad(H)\leq f(r)$. 	Thus for every graph $H$  whose $r$-subdivision $G$ is in $\mathscr C$ (thus in $\mathscr C_3$) we have $\ad(H)\leq f(r)$ hence \eqref{enum:OBE3} is satisfied.
\myqed\end{proof}

Strongly $\chi_p$-bounded classes of graphs can be also characterized by means of restricted homomorphism dualities (see \cite{Sparsity} for more background).
A {\em homomorphism} of a graph $\vec F$ to a graph $\vec G$ is a mapping $f:V(\vec F)\rightarrow V(\vec G)$ that preserve arcs: for every arc $uv$ of $\vec F$, $f(u)f(v)$ is an arc of $\vec G$. We denote by $\vec F\rightarrow\vec G$ the existence of a homomorphism of $\vec F$ to $\vec G$, and by $\vec F\nrightarrow\vec G$ the non-existence of such a homomorphism.
An {\em oriented graph} is a loopless directed graph with no circuits of length $2$, that is an orientation of an undirected graph. It is easily checked that the chromatic number of a graph $G$ is the minimum order of a loopless directed graph $\vec H$ (which can be required to be an oriented graph) such that some orientation $\vec G$ of $G$ satisfies $\vec G\rightarrow\vec H$. 

Let $\vec{\mathscr C}$ be a class of directed graphs. A directed graph $\vec{D}$ is a {\em restricted dual} of a directed graph $\vec{F}$ for the class $\vec{\mathscr C}$ if $\vec{F}\nrightarrow\vec{D}$ and, for every directed graph $\vec G\in\vec{\mathscr C}$ we have
	\begin{equation}
	\label{eq:rd}
		\vec F\nrightarrow\vec G\quad\iff\quad \vec G\rightarrow\vec D.	
	\end{equation}

The class $\vec{\mathscr C}$ has {\em all restricted dualities} if every directed connected graph $\vec{F}$ has a dual $\vec{D}$ for $\vec{\mathscr C}$.

\begin{lemma}
\label{lem:OBErd}
For a hereditary class $\mathscr C$ of graphs, let $\vec{\mathscr C}$ denote the class of all orientations of the graphs in $\mathscr C$. Then 
	the following are equivalent:
\begin{itemize}
	\item every connected acyclically oriented graph has a restricted dual for the class $\vec{\mathscr C}$;
	\item the class $\mathscr C$ has $\omega$-bounded expansion. 
\end{itemize}	
\end{lemma}
\begin{proof}
	Let $\vec{\mathscr C}$ be a hereditary class of oriented graphs closed under reorientation and let $\mathscr C$ be the underlying class of undirected graphs.
	
Assume $\mathscr C$ has $\omega$-bounded expansion.
Let $\vec F$ be a connected acyclically oriented graph. Let $t=2^{|\vec F|}$.
	As the class $\vec{\mathscr C_t}$ of all oriented graphs in $\vec{\mathscr C}$ with clique number at most $t$ has bounded expansion, it has all restricted dualities \cite{POMNIII}. Thus there exists $\vec D$ with $\vec F\nrightarrow\vec D$ such that for all $\vec G\in\vec{\mathscr C}_t$ the equivalence \eqref{eq:rd} holds.
	Let $\vec G\in \vec{\mathscr C}\setminus\vec{\mathscr C}_t$.
Then  $\omega(\vec{G})\geq 2^{|\vec F|}$ thus $\vec G$ contains a transitive tournament on $|\vec F|$ vertices. Hence $\vec F\rightarrow\vec G$. It follows that $\vec G\nrightarrow\vec D$, for otherwise we would deduce $\vec F\rightarrow\vec D$ by transitivity. It follows that $\vec D$ is a restricted dual of $\vec F$ for the class $\vec{\mathscr C}$.

	Assume for contradiction that every connected acyclically oriented graph $\vec F$ has a restricted dual $\vec D_F$ for the class $\vec{\mathscr C}$, but the class
	$\mathscr C$ does not have $\omega$-bounded expansion. We apply Lemma~\ref{lem:OBEi}.

	First assume that the class $\mathscr C$ is not $\chi$-bounded or that it includes all complete bipartite graphs. Then there is an integer $t$ such that  
	$\mathscr C$  contains graphs with arbitrarily large average degrees and clique number at most $t$.
	 Let $\vec T_{t+1}$ be the transitive tournament on $t+1$ vertices, and let $\vec{\mathscr C}'=\{\vec G\in\vec{\mathscr C}\mid \vec T_{t+1}\nrightarrow\vec G\}$.
	By assumption, $\vec{\mathscr C'}$ contains oriented graphs with 
	arbitrarily large average degree hence with arbitrarily large
	 oriented chromatic number, contradicting the property that 
	 every graph $\vec G\in\vec{\mathscr C}'$ satisfies $\vec G\rightarrow\vec D_{\vec T_{t+1}}$. 
	
	Otherwise, according to Lemma~\ref{lem:OBEi}, for some integer $q\geq 1$,  the class
	$\mathscr C$ contains the $q$-subdivisions of graphs with arbitrarily large average degree.
	Then, according to Lemma~\ref{lem:1subchr} the class
	$\mathscr C$ contains the $p$-subdivisions of graphs with arbitrarily large average degree, where $p=2q+1$.
	Let  $\mathscr D$ be the class of all  graphs, whose $p$-subdivision is in $\mathscr C$. As $p\geq 1$ and $\mathscr C$ is hereditary, the class $\mathscr D$ is monotone. By assumption, there are graphs in  $\mathscr D$ with arbitrarily large chromatic number. By \cite{Rodl1977}, $\mathscr D$ contains triangle free graphs with arbitrarily large chromatic number. Let $\mathscr D'$ be the class of  all  triangle free graphs in $\mathscr D$, and let $\vec{\mathscr C}'$ be the class of the $p$-subdivisions of all orientations of graphs in $\mathscr D'$.
		Let $\vec F$ be the $p$-subdivision of $\vec T_3$, and let $\vec D$ be its dual.
		Let $\vec{D}'$ be the directed graph with vertex set $V(\vec D)$, in which $uv$ is an arc if there exists in $\vec D$ a directed walk of length $p+1$ from $u$ to $v$. As $\vec F\nrightarrow\vec D$ the directed graph $\vec D'$ has no loops. Let $D'$ be the undirected graph underlying $\vec{D'}$.
		As $\vec F\nrightarrow\vec G$ for every $\vec G\in\vec{\mathscr C}'$ we deduce that every  $\vec G\in\vec{\mathscr C}'$ satisfies $\vec G\rightarrow\vec D$. It follows that for every $H\in\mathscr D'$ we have $H\rightarrow D'$, what contradicts the hypothesis that graphs in $\mathscr D'$ have arbitrarily large chromatic number.
\myqed\end{proof}
Note that every class $\mathscr C$ of undirected graphs with $\omega$-bounded expansion has all restricted dualities. 
Lemma~\ref{lem:OBErd} should be compared with the following characterization of bounded expansion classes.
\begin{lemma}
	For a class $\mathscr C$ of graphs, let $\vec{\mathscr C}$ denote the class of all orientations of the graphs in $\mathscr C$. Then
	the following are equivalent:
\begin{itemize}
	\item the class $\vec{\mathscr C}$ has all restricted dualities;
	\item the class $\mathscr C$ has bounded expansion. 
\end{itemize}	
\end{lemma}
\begin{proof}
	Assume $\mathscr C$ has bounded expansion. Then $\vec{\mathscr C}$ has all restricted dualities (see \cite{Sparsity}).
		
As $\mathscr C$ does not have bounded expansion there exists an integer $p$, such that for every integer $k$ there is a graph $H$ with chromatic number $k$, whose $p$-subdivision is a subgraph of a graph $G_H\in\mathscr C$ \cite{Sparsity}. Note that we can assume $p\geq 3$.
Let $\vec D_{\vec C_{p+1}}$ be a restricted dual of $\vec C_{p+1}$, the directed cycle of length $p+1$ and let $H$ be a graph with $\chi(H)>|\vec D_{\vec C_{p+1}}|$, whose $p$-subdivision is a subgraph of a graph $G_H\in\mathscr C$.
	Let $\vec H$ be an acyclic orientation of $H$ and let  
	$\vec{G}_H$ be an acyclic orientation of $G_H$ extending the orientation of the $p$-subdivision of $H$ inherited from $\vec H$.
	 Let $\vec D$ be the directed graph with vertex set $V(\vec D_{\vec C_{p+1}})$, in which $uv$ is an arc if there exists in $\vec D_{\vec C_{p+1}}$ a directed walk of length $p+1$ from $u$ to $v$. As $\vec C_{p+1}\nrightarrow \vec D_{\vec C_{p+1}}$, the directed graph $\vec D$ is loopless. As $\vec C_{p+1}\nrightarrow \vec{G}_H$ we have $\vec{G}_H\rightarrow  \vec D_{\vec C_{p+1}}$ thus $\vec H\rightarrow\vec D$, which implies $\chi(H)\leq |\vec D|$.
\myqed\end{proof}

\section{Weakly $\chi_p$-bounded classes of graphs}
\label{sec:weak}
We illustrate the difference between the notions of strongly $\chi_p$-bounded classes and weakly $\chi_p$-bounded classes with the following example.

\begin{proposition}
\label{prop:biw}	
	The class $\mathscr B$ of all complete bipartite graphs is weakly $\chi_p$-bounded but not strongly $\chi_p$-bounded.
\end{proposition}
\begin{proof}
	The class $\mathscr B$ is clearly not strongly $\chi_p$-bounded as it is triangle free but has unbounded $\chi_2$.
	However, the class $\mathscr B$ is obviously weakly $\chi_1$-bounded (i.e. $\chi$-bounded).
	Let $p\geq 2$.
	Let $K_{s,t}$ be a complete bipartite graphs with $s\leq t$. Then $\chi_p(K_{s,t})\leq {\rm td}(K_{s,t})\leq s+1$. Moreover, $\max_{H\in {\rm TM}_{p-1}(G)}\omega(H)\geq \max_{H\in {\rm TM}_{1}(G)}\omega(H)\geq \sqrt{s}$. It follows that for every integer $p\geq 2$ and every complete bipartite graph $K_{s,t}$ we have
\[
	\chi_p(K_{s,t})\leq\max_{H\in {\rm TM}_{p-1}(K_{s,t})}\omega(H)^2,
\]
thus $\mathscr B$ is weakly $\chi_p$-bounded.
\myqed\end{proof}

The class of all $1$-subdivision of graphs is an example of a $\chi$-bounded class (as it includes only bipartite graphs) that is not weakly $\chi_p$-bounded (as the class is $C_4$-free, while $\chi_s$ is unbounded by Lemma~\ref{lem:chis-sub}). This suggests that  $\chi_p(G)$ should be related to the chromatic number of shallow topological minors of $G$, which is the subject of the next two results.
\begin{lemma}
\label{lem:chiptm}
	Let $G$ be a graph and let $p$ be a positive integer. Then
	\begin{equation}
		\chi_p(G)\geq\chi({\rm TM}_{p-1}(G))^{\frac{1}{p}}.
	\end{equation}
\end{lemma}
\begin{proof}
	For $p=1$ there is nothing to be proved. 
	
	For $p=2$, the proof of Lemma~\ref{lem:chis-sub} can be easily modified to give the result:
	Let $G$ be a graph and let $H\in {\rm TM}_{1}(G)$.
	We shall consider $V(H)$ as a subset of $V(G)$.
	Let $c$ be a star-coloring (i.e. a $\chi_2$-coloring) of $G$ with $k=\chi_2(G)$ colors. 	 
	 Let $H_0\subseteq H$ be the spanning subgraph of $H$ with edge set $E(H_0)=\{uv\in E(H): c(u)=c(v)\}$.
	  Then every connected component of $H_0$ is monochromatic under $c$ and all the edges of $H_0$ correspond to paths of length $2$ in $G$ (i.e. are $1$-subdivided in $G$).
	According to \cite[Lemma 4]{wood2005acyclic},  
	 the minimum number of colors in a star colouring of the $1$-subdivision $H_0'$ of $H_0$ in which the original vertices are monochromatic is $\chi'(H_0)+1$, where $\chi'(H_0)$ denotes the edge chromatic number of $H_0$. Hence the edge $\chi'(H_0)\leq k-1$. Thus $\Delta(H_0)\leq k-1$ and $\chi(H_0)\leq k$ by Brooks theorem. Let $\phi$ be a vertex $k$-coloring of $H_0$. Now color each vertex $v$ of $H$ by the pair $(c(v),\phi(v))$. Consider an edge $uv\in E(H)$. If $uv\in E(H_0)$ then $\phi(u)\neq\phi(v)$. If $uv\notin E(H_0)$ then $c(u)\neq c(v)$. Thus we have a $k^2$-coloring of $H$, and $\chi(H)\leq \chi_2(G)^2$.


	Now assume $p\geq 3$.
	Let $G$ be a graph and let $H\in{\rm TM}_{p-1}(G)$.
	The vertices of $H$ naturally correspond to vertices $a_1,\dots,a_{|H|}$ of $G$, and to each edge $a_ia_j$ of $H$ corresponds a path $P_{i,j}$ of $G$ with length at most $p$ linking $a_i$ and $a_j$.
		Consider a $\chi_p$-coloring of $G$ with $k=\chi_p(G)$ colors (taken in $[k]$). This coloring naturally defines a coloring of the vertices of $H$. 	
		By pigeon-hole principle, there exists a color $c\in [k]$ such that the subgraph $H_c$ of $H$ induced by vertices colored $c$ has chromatic number at least $\lceil\chi(H)/k\rceil$. 
		It follows that $H_c$ contains an induced subgraph $\hat H_c$ with 
		average degree at least $\lceil\chi(H)/k\rceil-1$.
		 For each edge $e=a_ia_j$ of $\hat H_c$ 
		  we denote by $\gamma(e)$ the set of all the colors present (in $G$) on the path $P_{i,j}$. To each $(p-1)$-subset $I$ of $[k]\setminus\{c\}$ corresponds a subset $E_I$ of edges $e$ of $\hat H_c$ with $\gamma(e)\subseteq I\cup\{c\}$. By pigeon-hole principle, there exists a subset $I$ such that the 
		 		 average degree of the subgraph $H_{c,I}$ of $\hat H_c$ induced by $E_I$ is at least $(\lceil\chi(H)/k\rceil-1)/\binom{k-1}{p-1}$. The graph $H_{c,I}$ defines a subgraph $G_{c,I}$ of $G$ by taking the union of all the paths $P_{i,j}$ for $a_ia_j\in E(H_{c,I})$. By construction, the vertices of $G_{c,I}$ are colored by colors in $I\cup\{c\}$, which is a subset of $p$ colors. Thus ${\rm td}(H_{c,I})\leq{\rm td}(G_{c,I})\leq p$. It follows that $H_{c,I}$ has average degree less than $2p-2$.
	 		 It follows that 
\begin{align*}
	2p-2&>\frac{\lceil\chi(H)/k\rceil-1}{\binom{k-1}{p-1}}.
\intertext{Thus}
	\chi(H)&<\left(2(p-1)\binom{k-1}{p-1}+1\right)k<2p^2\binom{k}{p}\leq \frac{2p}{(p-1)!}k^p
\end{align*}
In particular, if $p>3$ then $\chi(H)<\chi_p(G)^p$.

So assume $p=3$, and let $C$ be a connected component of $G_{c,I}$ with maximal average degree. If we remove the root $r$ of $C$ we are left with a star forest. If  $r$ is not colored $c$ it
follows that $r$ has degree at most $2$ hence $C$ contains a most one cycle thus the average degree of $C$ is at most $2$.
So assume that $r$ is colored $c$. Assume some connected component of $H_{c,I}-r$ contains two adjacent vertices $u$ and $v$, at least one of them (say $u$) being adjacent to $r$. Then in a connected component $C-r$ we have a path of length $4$ (at least a subdivision vertex for the edge $ru$, the vertex $u$, at least a subdivision vertex for the edge $uv$, then the vertex $v$), contradicting the hypothesis that the connected components of $C-r$ are stars. Hence $H_{c,I}$ is a star, 
and its average degree is at most $2$. Thus $\chi(H)<(2\binom{k-1}{2}+1)k<k^3$. 
	\myqed\end{proof}

We deduce the following generalization of Lemma~\ref{lem:chis-sub}, which is of independent interest.
\begin{theorem}
\label{thm:chip-sub}
	Let $p$ be a non negative integer.
	Let $G$ be a graph and let $G^{(p)}$ be its $p$-subdivision. Then
\begin{equation}
	\chi(G)^{\frac{1}{p+1}}\leq \chi_{p+1}(G^{(p)})\leq \max(\chi(G),p+2).
\end{equation}
\end{theorem}
\begin{proof}
	If $p=0$ the statement obviously holds.
	According to Lemma~\ref{lem:chiptm} we only have to prove the inequality $\chi_{p+1}(G^{(p)})\leq \max(\chi(G),p+2)$. Consider a proper coloring of $G$ with $\chi(G)$ colors and transfer the colors on the corresponding vertices of $G^{(p)}$. Then, for each edge $uv$ of $G$,  we color the $p$ internal vertices of the  path $P_{uv}$ of $G^{(p)}$ corresponding to the edge $uv$ of $G$ by distinct $p$ colors that are also distinct from the color of $u$ and the color of $v$. It is easily checked that every subset of $k\leq p+1$ colors then induce a subgraph of $G^{(p)}$ with tree-depth at most $k$. 
\myqed\end{proof}

We now state two lemmas that will lead to the proof of Theorem~\ref{thm:weak}.

\begin{lemma}
\label{lem:1}
	Let $\mathscr C$ be a hereditary class of graph and let $p$ be a positive integer.
	If $\chi_p(G)\leq f(\omega({\rm TM}_{p-1}(G)))$ holds for all $G\in\mathscr C$ then 
	${\rm ITM}_{p-1}^e(\mathscr C)$ is $\chi$-bounded. 
\end{lemma}
\begin{proof}
	Without loss of generality, we can assume that $f$ is non-decreasing. 
	Let $H\in {\rm ITM}_{p-1}^e(\mathscr C)$. Then there exists $G\in\mathscr C$ such that $H^{(p-1)}\subseteq_i G$. As $\mathscr C$ is hereditary we deduce $H^{(p-1)}\in\mathscr C$. 
	
	Hence we have
\begin{align*}
	\chi(H)&\leq \chi({\rm TM}_{p-1}(H^{(p-1)}))&\text{(as $H\in {\rm TM}_{p-1}(H^{(p-1)})$)}\\
	&\leq \chi_p(H^{(p-1)})^p&\text{(by Lemma~\ref{lem:chiptm})}\\
	&\leq f(\omega({\rm TM}_{p-1}(H^{(p-1)})))^p&\text{(by Lemma assumption)}\\
	&\leq f(\omega(H))^p&\text{(as $\omega({\rm TM}_{p-1}(H^{(p-1)}))=\omega(H)$)}
\end{align*}
It follows that ${\rm ITM}_{p-1}^e(\mathscr C)$ is $\chi$-bounded.
	\myqed\end{proof}

The following is an easy but useful lemma.
\begin{lemma}
	Let $\mathscr D$ be a monotone class of graph.
	If $\mathscr D$ is $\chi$-bounded then
	$\chi(\mathscr D)<\infty$. 
\end{lemma}
\begin{proof}
Assume for contradiction $\chi(\mathscr D)=\infty$.
By \cite{Rodl1977} it follows  that $\mathscr D$ contains triangle-free graphs with arbitrarily large chromatic number, contradicting the assumption that $\mathscr D$ is $\chi$-bounded.\myqed\end{proof}

As ${\rm ITM}_{p}^e(\mathscr C)$ is obviously monotone for $p\geq 1$ we deduce the following strengthening of Lemma~\ref{lem:1}.
\begin{corollary}
\label{cor:1}
	Let $\mathscr C$ be a hereditary class of graphs and let $p$ be a positive integer.
	If $\chi_{p+1}(G)\leq f(\omega({\rm TM}_{p}(G))$ holds for all $G\in\mathscr C$ then 
	$\chi({\rm ITM}_{p}^e(\mathscr C))<\infty$. 
\end{corollary}

\begin{potweak}
	\eqref{enum:w1}$\Rightarrow$\eqref{enum:w2}: The $\chi$-boundedness of $\mathscr C$ is simply the case $p=1$ of \eqref{enum:w1} and the other cases directly follow from Corollary~\ref{cor:1}.
	
	\eqref{enum:w2}$\Rightarrow$\eqref{enum:w3}:  Assume for contradiction that ${\rm ITM}_{r}^e(\mathscr C)$ has unbounded average degree. According to Lemma~\ref{lem:1subchr} the chromatic number of graphs in ${\rm ITM}_{2r+1}^e(\mathscr C)$ is unbounded, contradicting \eqref{enum:w2}.
	
	\eqref{enum:w3}$\Rightarrow$\eqref{enum:w4}: Let $s$ be a positive integer and let $\mathscr C_s=\{G\in\mathscr C\mid\bomega(G)\leq s\}$. As $\ad({\rm ITM}_1^e(\mathscr C))<\infty$ we deduce from Lemma~\ref{lem:dense1sub} that $\ad(\mathscr C)<\infty$.
	Then, according to Lemma~\ref{lem:zd} we deduce $\ad({\rm TM}_k(\mathscr C))<\infty$ for all $k$. It follows that $\mathscr C_s$ has bounded expansion (by Lemma~\ref{lem:BEtm}) hence $\chi_p(\mathscr C)<\infty)$ for all integers $p$ (by Lemma~\ref{lem:BEchip}).
	As this holds for each integer $s$ we deduce that there exists a function $f_p$ for each integer $p$ such that $\chi_p(G)\leq f_p(\bomega(G))$ holds for every $G\in\mathscr C$.
	
	\eqref{enum:w4}$\Rightarrow$\eqref{enum:w5}: For each integer $t$, the class $\{G\in\mathscr C\mid K_{s,s}\not\subseteq_i G\text{ and }\omega(G)\leq t\}$ is included (by Ramsey theorem) in the class $\{G\in\mathscr C\mid K_{s^t,s^t}\not\subseteq G\}$, which has bounded expansion. It follows that the class $\{G\in\mathscr C\mid K_{s,s}\not\subseteq_i G\}$ is strongly $\chi_p$-bounded by Theorem~\ref{thm:strong}.
	
	\eqref{enum:w5}$\Rightarrow$\eqref{enum:w1}: assume $\mathscr C$ is $\chi$-bounded (i.e. there is a function $f$ with $\chi(G)\leq f(\omega(G))$ for all $G\in\mathscr C$) and that for each positive integer $s$ the class $\{G\in\mathscr C\mid K_{s,s}\not\subseteq_i G\}$ is strongly $\chi_p$-bounded (i.e. there are functions $g_{p,s}$ with  $\chi_p(G)\leq g_{p,s}(\omega(G))$ for all graphs $G\in\mathscr C$ with $K_{s,s}\not\subseteq_i G$).
		Let $G\in\mathscr C$. We have $\chi_1(G)\leq f(\omega(G))$ by assumption. Let $p>1$ and let $k=\omega({\rm TM}_{p-1}(G))$. Note that $G$ has clique number at most $k$ and that it does not contain any induced $K_{k+1,k+1}$ (as $K_{k+1}\in {\rm TM}_{p-1}(K_{k+1,k+1})$).
		Thus $\chi_p(G)\leq g_{p,k+1}(k)$.
		We conclude that $\mathscr C$ is weakly $\chi_p$-bounded.
		\qed\end{potweak}
	
With Theorem~\ref{thm:weak} at hand we can give some examples of weakly $\chi_p$-bounded classes.
First we note the following
 direct consequence of Lemma~\ref{lem:dvo-sub}.

\begin{itemize}
	\item classes with no holes of length greater than $\ell$ are weakly $\chi_p$-bounded,
		as they are $\chi$-bounded \cite{chudnovsky2017induced} and they exclude all subdivisions of $C_{\ell+1}$;
	\item classes excluding all subdivisions of some tree $T$ are weakly $\chi_p$-bounded,  as these classes are $\chi$-bounded \cite{Scott1997};
	\item The class of graphs with no holes of length equal to $0\bmod  \ell$  is weakly $\chi_p$-bounded. Indeed this class is $\chi$-bounded \cite{scott2019induced} and, for each integer $r$, the class ${\rm ITM}_r^e(\mathscr C)$ contains no cycle of length $0\bmod \ell$ hence is also $\chi$-bounded. 
\end{itemize}

We now give some further examples, which show a surprising robusteness of the notion of weak $\chi_p$-boundedness.
For this, we shall need the following result (stated as Lemma~\ref{lem:SBEchi}), which is a direct corollary of the main theorem of \cite{SBE_TOCL} (see   \cite{RW_SODA_arxiv} for a formal derivation). A {\em first-order transduction} $\mathsf T$ is a pair $(\eta(x,y), \nu(x))$ of first-order formulas in the language of vertex-colored graphs. A class $\mathscr D$ is a {\em first-order transduction} of a class $\mathscr C$ if there exists a first-order transduction $\mathsf T=(\eta(x,y), \nu(x))$, where $\eta$ is symmetric and for every graph $G\in\mathscr D$ there exists a graph $H\in\mathscr C$ and a vertex-coloring $H^+$ of $H$ such that 
\begin{itemize}
	\item $V(G)$ is the set $\nu(H^+)$ of vertices of $H$ that satisfy $\nu$ in $H^+$;
	\item $E(G)$ is the set $\eta(H^+)\cap \nu(H^+)^2$ of pairs vertices of $G$ that satisfy $\nu$ in $H^+$.
\end{itemize}
As an example, consider the class of all map graphs.
Recall that a \emph{map graph} is the intersection graphs of finitely many simply connected and internally disjoint regions of the plane, and that it can be obtained as induced subgraph of the square of a bipartite planar graph. The class of map graphs is a first-order transduction of the class of planar graphs, as witnessed by the transduction $\mathsf T=(\eta,\nu)$ where we consider a black/white coloring of the vertices, $\nu(x):={\rm Black}(x)$, and 
	$\eta(x,y):=(\exists z)\ E(x,z)\wedge E(z,y)$.
	
	A class $\mathscr C$ has {\em structurally bounded expansion} if it is a first-order transduction of a class with bounded expansion. For instance, the class of maps, which is a first-order transduction of the class of planar graphs, has structurally bounded expansion.

\begin{lemma}[\cite{SBE_TOCL,RW_SODA_arxiv}]
\label{lem:SBEchi}
	Every structurally bounded expansion class is linearly $\chi$-bounded.
\end{lemma}

We complement this by
\begin{proposition}
\label{prop:SBEchip}
	Every structurally bounded expansion class is weakly $\chi_p$-bounded for every $p$.
\end{proposition}	
\begin{proof}
	Let $\mathscr C$ be a structurally bounded expansion class.
	According to Lemma~\ref{lem:SBEchi} $\mathscr C$ is linearly $\chi$-bounded.
	Moreover, the classes ${\rm ITM}_{r}^e(\mathscr C)$ are  transductions of $\mathscr C$ hence are classes with structurally bounded expansion thus are also linearly $\chi$-bounded. By Theorem~\ref{thm:weak}-\eqref{enum:w2} we deduce that $\mathscr C$ is weakly $\chi_p$-bounded.
\myqed\end{proof}

For the next example we recall the notion of low $\sf P$-covers  \cite{SBE_TOCL}.
Let $\sf P$ be a hereditary class property.
A class $\mathscr C$ has {\em low $\sf P$-covers} if, for every positive integer $p$ there exists a class $\mathscr D_p$ with property $\mathsf P$ and an integer $n_p$ such that for every graph $G\in\mathscr C$ there exists family $\mathcal F$ of at most $n_p$ subsets of vertices of $G$ with the following property:
\begin{enumerate}
	\item every subset $X$ of at most $p$ vertices of $G$ is included in some set in $\mathcal F$;
	\item every subgraph of $G$ induced by a set in $\mathcal F$ belongs to $\mathscr D_p$.
\end{enumerate}

If $f$ is a graph invariant, a class $\mathscr C$ has {\em low $f$-covers}  if it has low $\mathsf P_f$-covers for the property $\mathsf P_f$ expressing that $f$ is bounded on the class (i.e. the corresponding classes $\mathscr D_p$ can be chosen of the form $\{G: f(G)\leq C_p\}$ for some constant $C_p$ depending on $p$).

A class $\mathscr C$ has {\em \bomega-bounded expansion} if there exists a function $f$ such that for every integer $r$ and every graph $G\in\mathscr C$ we have
\[
	\ad({\rm TM}_t(G))\leq f(\bomega(G),r).
\]
Note that it follows from Theorem~\ref{thm:weak} that a class is weakly $\chi_p$-bounded if and only if it is $\chi$-bounded and it has \bomega-bounded expansion.

The {\em twin-width} invariant has been recently introduced \cite{bonnet2020twinI, bonnet2020twinII, bonnet2020twinIII,
TWW4_arxiv, TWWP-arxiv}.
Classes with bounded twin-width include
proper minor-closed classes and bounded rank-width graphs, and the property of having bounded twin-width is preserved by first-order transductions.

\begin{proposition}
	Every class with low twin-width covers is weakly $\chi_p$-bounded.
\end{proposition}
\begin{proof}
	Let $\mathscr C$ be a class with low twin-width covers. 
	Then there exist integers $n_1,t_1$ such that the vertex set of every graph $G$ in $\mathscr C$ can be partitioned into at most $n_1$ parts $V_1,\dots,V_{n_1}$, each inducing a subgraph with twin-width at most $t_1$. As classes with bounded twin-width  are $\chi$-bounded \cite{bonnet2020twinIII}, there is a function $f$ such that $\chi(G[V_i])\leq f(\omega(G[V_i])$ for each $1\leq i\leq n_1$. Hence $\chi(G)\leq \sum_{i=1}^{n_1}\chi(G[V_i])\leq n_1\,f(\omega(G))$. Thus $\mathscr C$ is $\chi$-bounded.
	As a class with bounded twin-width and bounded {\bomega} has bounded expansion  \cite{bonnet2020twinII}, classes with 
	low twin-width covers and bounded {\bomega} have low bounded expansion covers hence have bounded expansion \cite{POMNI}.
\myqed\end{proof}

\begin{proposition}
\label{prop:itmclosed}
	A class $\mathscr C$ closed under induced topological minors is weakly $\chi_p$-bounded if and only if it is $\chi$-bounded.
\end{proposition}
\begin{proof}
	Let $\mathscr C$ be a class of graphs closed under induced topological minors.
	
	If $\mathscr C$ is weakly $\chi_p$-bounded then it is  $\chi$-bounded.
	
	Conversely, assume $\mathscr C$ is $\chi$-bounded.
	For every positive integer $r$, the class ${\rm ITM}_r^e(\mathscr C)$ is included in $\mathscr C$. Hence the classes ${\rm ITM}_r^e(\mathscr C)$ are $\chi$-bounded. By Theorem~\ref{thm:weak} it follows that $\mathscr C$ is weakly $\chi_p$-bounded.
\myqed\end{proof}

Very recently, Davies~\cite{davies2020} announced that proper vertex-minor-closed classes are $\chi$-bounded. 
More generally, we now show that it follows that such classes are weakly $\chi_p$-bounded.

Recall that the {\em local complementation} at a vertex $v$ of a graph $G$ is the operation of replacing the subgraph induced by the neighborhood of $v$ by its complement, and that the resulting graph is denoted by $G\ast v$. A graph $H$ is a {\em vertex-minor} of a graph $G$ if it can be derived from $G$ by applying a sequence of local complementations and vertex deletions. A class is {\em vertex-minor-closed} if every vertex-minor of a graph in the class also belongs to the class; it is {\em proper} if it does not include all graphs.

	It is easily checked that if a subdivision of a graph $H$ is an induced subgraph of a graph $G$ then $H$ is a vertex-minor of $G$. It follows that  proper vertex-minor closed classes of graphs are closed under induced topological minors. It follows then from  Proposition~\ref{prop:itmclosed} that $\chi$-boundedness of proper vertex-minor closed classes of graphs imply weak $\chi_p$-boundedness.

On another hand, it is interesting to compare the notions of weakly $\chi_p$-bounded class and nowhere dense class. 

\begin{proposition}
	Let $\mathscr C$ be a weakly $\chi_p$-bounded class of graphs. Then $\mathscr C$ is nowhere dense if and only if $\mathscr C$ has bounded expansion.
\end{proposition}
\begin{proof}
	If $\mathscr C$ has bounded expansion it is obviously nowhere dense. So assume $\mathscr C$ is a nowhere dense weakly $\chi_p$-bounded class.
	As $\mathscr C$ is nowhere dense, we have $\omega({\rm TM}_r(\mathscr C))<\infty$ for every non-negative integer $r$, thus $\chi_p$ is bounded on $\mathscr C$ for each integer $p$ hence $\mathscr C$ has bounded expansion.
\myqed\end{proof}

\section{Examples}
\label{sec:appli}
In this section, we provide some examples to motivate the introduction of the notion of classes with $\omega$-bounded expansion.

\subsection{Fixed parameter tractability}
The first application builds on the fixed parameter linear time algorithm for first-order model checking on bounded expansion classes proposed by Dvo\v r\'ak, Kra{\v l}, and Thomas~\cite{DKT,DKT2}. It is an immediate consequence of this result that any first order-sentence $\varphi$ can be tested in graphs $G$ in a class $\mathscr C$ with $\omega$-bounded expansion in time
	$f(|\varphi|,\omega(G))\,|G|$, thus first-order model checking is FPT on $\mathscr C$ when parametrized by both the length of the sentence and the clique-number of the graph. We now prove that the dependence to the clique-number can be avoided if we restrict the sentences expressing a property preserved when taking a supergraph, like existential positive sentences. 
 
\begin{proposition}
	Let $\mathscr C$ be a class with $\omega$-bounded expansion.
	Then for every first-order sentence $\phi$ expressing a monotone  property (in the sense that every supergraph of a model of $\phi$ is a model of $\phi$) there is a linear time algorithm $\mathsf A$ that checks whether $G\in\mathscr C$ satisfies $\phi$. Moreover, if $G\in\mathscr C$ satisfies 
	$\phi$ then the output of the algorithm $\mathsf A$ is a minimal subset $X$ with $|X|\leq t$ and $G[X]\models \phi$.
\end{proposition}
\begin{proof}
	As $\mathscr C$ has $\omega$-bounded expansion there exists a function $f$ such that $\ad(G)\leq f(\omega(G))$ for every $G\in\mathscr C$.
	Let $\phi$ be a first-order sentence expressing a monotone property, in the sense that if $G_1$ is a subgraph of a graph $G_2$ and $G_1\models \phi$ then $G_2\models \phi$. 
	If $\phi$ has no finite model\footnote{Note that this is equivalent to the property that no finite complete graph satisfies $\phi$, which is decidable.}, then the algorithm outputs `No' for every input instance. Otherwise, let $t$ be the minimum integer such that  $G\models K_t$. 
	For $G\in\mathscr C$ we first compute a topological ordering of $G$, from which we deduce the degeneracy of $G$. If the degeneracy of $G$ is smaller than $f(t)$ then $\omega(G)\leq f(t)$. Thus $G$ belongs to $\mathscr C_{f(t)}$, which is a class with bounded expansion in which $\phi$ can be tested in linear time~\cite{DKT, DKT2}. We  do so easily using a low tree-depth decomposition with parameter $t$, which provides us a witness. Otherwise, $\omega(G)>t$ thus $G\models\phi$.
	By using the topological order, we find an induced subgraph $H$ of $G$ that is $f(t)$ degenerate. Then a clique of size $t$ can be found in $H$ in time $\binom{f(t)}{t}|H|$.
\myqed\end{proof}

We leave the following as a problem.  
\begin{problem}
	Is first-order model checking is FPT on hereditary classes with $\omega$-bounded expansion?
\end{problem}

The second application provides us an effective strengthening of restricted dualities.
Recall that a \emph{restricted dual} of a graph $F$ with respect to a class of graphs $\mathscr C$ is a graph $D_F$  that is not homomorphic to $F$ (i.e. $F\nrightarrow D_F$) and such that, for every graph $G\in\mathscr C$, we have
\[
	F\nrightarrow G\qquad\iff\qquad G\rightarrow D_F.
\]
Note that for every graph $G$ (including graphs out of $\mathscr C$) we have $(F\rightarrow G)\Longrightarrow (G\nrightarrow D_F)$, for otherwise we would have $F\rightarrow G\rightarrow D_F$.
It was proved in \cite{POMNIII} that if $\mathscr C$ is a class with bounded expansion then every connected graph $F$ has a restricted dual with respect to $\mathscr C$. This property is the core of the next proposition.

\begin{proposition}
Let $\mathscr C$ be a class with $\omega$-bounded expansion and unbounded clique number.

	Then there exists infinitely many $H$-coloring problems that are not equivalent on $\mathscr C$ and that can be solved in linear time on $\mathscr C$ with witness, in the following sense: 
	
	for each of the $H$-coloring problems, there exists a constant $C$ and a linear time algorithm $\mathsf A$ such that for input graph $G\in\mathscr C$ the algorithm $\mathsf A$ outputs 
	\begin{itemize}
		\item either a homomorphism $f:G\rightarrow H$,
		\item or a subset $X$ of at most $C$ vertices of $G$ with $G[X]\nrightarrow H$.
	\end{itemize}
	
\end{proposition}
\begin{proof}
	As $\mathscr C$ has $\omega$-bounded expansion there exists a function $f$ such that every $G\in\mathscr C$ is $f(\omega(G))$-degenerate. 
	For an integer $k$, let $\mathscr C_k$ denote the subclass of $f(k)$-degenerate graphs in $\mathscr C$.  

	Let $F$ be an arbitrary connected graph, let $C=|F|$, and let $H$ be the dual of $F$ with respect to $\mathscr C_{C}$ constructed in \cite{POMNIII}. Let us prove that the $H$-coloring problem can be solved in linear time on $\mathscr C$ with witness.
	We first perform a topological sort on $G$ (in linear time) and determine whether $G$ is $f(C)$-degenerate. If not, we can easily  extract in linear time (by an easy modification of the topological sort algorithm) an induced subgraph $G'$ of $G$ that is $f(C)$-degenerate but not $f(C-1)$-degenerate. Thus $\omega(G')>C-1$, which implies $F\rightarrow G'$ thus $G'\nrightarrow H$ (as $G'\in\mathscr C_C$). A subset $X$ with $|X|=C$ and $F\subseteq G'[X]=G[X]$ (hence $G[X]\nrightarrow H$) is extracted in linear time \cite{POMNII} and output.
	Otherwise, the graph $G$ belongs to the bounded expansion class $\mathscr C$.	Then we can check $F\rightarrow G$ in linear time \cite{POMNII}. If $F\rightarrow G$, the algorithm outputs a subset $X$ with $|X|=C$ and $F\subseteq G[X]$ (hence $G[X]\nrightarrow H$).
	Otherwise,  $F\nrightarrow G$. Then $\omega(G)\leq |F|$ and 
	we can use a low tree-depth decomposition with parameter $|F|$ to compute a coloring $G\rightarrow H$. (This follows from the construction in  \cite{POMNIII}.)

 As $\mathscr C$ has unbounded clique numbers, there exists infinitely many graphs $G_i\in\mathscr C$ ($i\in\mathbb N$) with strictly increasing clique numbers. Let $H_i$ be the restricted dual of $K_{\omega(G_i)}$ with respect to $\mathscr C_{\omega(G_i)}$. For every $j<i$ we have $G_j\in\mathscr C_{\omega(G_i)}$ and $K_{\omega(G_i)}\nrightarrow G_j$ thus $G_j\rightarrow H_i$. If $j\geq i$ then $\omega(G_j)\geq \omega(G_i)$ thus $G_j$ contains a clique of size $\omega(G_i)$ hence $G_i\nrightarrow H_i$. Altogether, $G_j\rightarrow H_i$ if and only if $j<i$. This witnesses that the $H_i$-coloring problems are not equivalent on $\mathscr C$.
\myqed\end{proof}

\begin{remark}
	Triangle-free even hole-free graphs are $2$-degenerate \cite{chudnovsky2019even} hence are $3$-colorable. This can be extended using restricted dualities:
	
	For every odd integer $g\geq 5$ there exists a $3$-colorable graph $D_g$ with odd girth $g$, such that for every even hole-free graph $G$ we have:
\[
	\text{\rm odd-girth}(G)\geq g\iff G\rightarrow D_g.
\]
\end{remark}

\subsection{Holes in even-hole free graphs}
As a further example of application of the notion of strongly $\chi_p$-bounded class, we consider the number $h_g(G)$ of holes of length $g$ in an even hole-free graph $G$. We shall make use of the following technical lemma.

\begin{lemma}[\cite{Taxi_tdepth}]
\label{lem:td_tech}
	There exists a function $F:\mathbb N\rightarrow\mathbb N$ with the following property: for every graph $G$ with tree-depth at most $t$ and at least $F(t)$ vertices there exists a partition $(A,X,B)$ of $V(G)$ such that:
	\begin{itemize}
		\item there is no edge between vertices in $A$ and vertices and $B$;
		\item there is a bijection $f:A\rightarrow B$ inducing an isomorphism of $G[A\cup X]$ and $G[B\cup X]$.
	\end{itemize}
\end{lemma}

\begin{lemma}
\label{lem:htd}
	For every odd integer $g$ there exists a constant $C_1(g)$ such that every even hole-free graph $G$ with tree-depth at most $g$ contains at most $C_1(g)\,|G|$ holes with length $g$. 
\end{lemma}
\begin{proof}
	Let $F$ be a function fulfilling the requirements of Lemma~\ref{lem:td_tech}.
	We construct three sequences $G_0,G_1,\dots,G_{k-1}$, $K_0,\dots,K_{k-1}$ and $H_1,\dots,H_k$ of induced subgraphs of $G$ inductively as follows: we let $G_0=G$ and, as long as $|G_i|>F(g)$ (for $i\geq 0$) there exists, according to Lemma~\ref{lem:td_tech}, a partition $(A,X,B)$ of $V(G_i)$ such that 
	there is no edge between $A$ and $B$, and there is a bijection $f:A\rightarrow B$ inducing an isomorphism of $G_i[A\cup X]$ and $G_i[B\cup X]$. We further consider a triple $(A,X,B)$ with $A$ minimal (for inclusion).
	Let $K$ be the subset of vertices of $X$ with at least one neighbour in $A$. By minimality of $A$, the subgraph $G_i[K\cup A]$ is connected. 
	Assume $u,v$ are distinct non-adjacent vertices of $K$ and let $P_A$ be a shortest path linking $u$ and $v$ in $G_i[K\cup A]$. Then $V(P_A)\cup f(V(P_A))$ induce an even hole of $G_i$, contradicting the assumption that $G$ (thus $G_i$) is even hole-free. Thus $K$ is a clique.
	
	Let $K_i=K$, $H_{i+1}=G_i[K\cup A]$ and  $G_{i+1}=G_i[X\cup B]$. 
	
	If $|G_i|\leq F(g)$ we stop the process (and we let $H_{i+1}=G_i$, and $k=i+1$).
	Note that $V(G)=\bigcup_{i\leq k}V(H_i\setminus K_{i-1})$.
	Moreover, as holes cannot cross clique separators, every hole of $G$ is fully included in one $G_i$. Thus
	\[
		h_g(G)=\sum_{i=1}^k h_g(G_i).
	\]
	Moreover, $k\leq |G_i|$ as each iteration removes at least one vertex.
	Let $C_1(g)$ be the maximum number of $h_g(G)$ for even hole-free graphs $G$ with tree-depth at most $g$ and at most $F(g)$ vertices.
	Then we have $h_g(G)\leq C_1(g)\,|G|$.
		\myqed
\end{proof}

\begin{theorem}
\label{thm:holes}
	The ratio
	$h_g(G)/|G|$ is bounded by a function of $g$ and $\omega(G)$ on even hole-free graphs. 
\end{theorem}

\begin{proof}
	First note that we only have to consider odd integers $g>3$ as otherwise $h_g(G)=0$.
	As the class of even hole-free graphs is strongly $\chi_p$-bounded there exists for every integers $g,\omega$ a constant $N_{g,\omega}$ such that every even hole-free graph $G$ with $\omega(G)=\omega$ has a vertex coloring $c:V(G)\rightarrow [N_{g,\omega}]$ such that every subset $I$ of $g$ colors induce a subgraph $G_I$ with tree-depth at most $g$. Obviously, the number of holes with length $g$ in $G$ is at most the sum of the number of holes with length $g$ in the subgraphs $G_I$ when $I$ ranges over all subsets of size $g$ of $[N_{g,\omega}]$.
	Thus, by Lemma~\ref{lem:htd} we have
\[	h_g(G)=\sum_{I\in\binom{[N_{g,\omega}]}{g}}h_g(G_I)
		\leq \sum_{I\in\binom{[N_{g,\omega}]}{g}}C_1(g)|G_I|
		\leq 		\binom{N_{g,\omega}}{g}\,C_1(g)\,|G|.
\]
\myqed
\end{proof}
Note that in Theorem~\ref{thm:holes} the dependence in $g$ and $\omega$ is necessary as, for even $\omega$ and odd $g>3$, the disjoint union $G$ of $(2n/g\omega)$ copies of $C_g[K_{\omega/2}]$ is even hole-free and
\[
	h_g(G)=\frac{1}{g}\left(\frac{\omega}{2}\right)^{g-1}\,|G|.
\]

\section{Conclusion}
In this paper we investigate the $\chi_p$ boundedness in terms of cliques in the graph (i.e. cliques at depth $0$ in the case of strong $\chi_p$-boundedness) and in terms of cliques of shallow minors at depth $p-1$
(in the case of weak $\chi_p$-boundedness). One may wonder whether one could investigate $\chi_p$  boundedness in terms of cliques in the shallow minors at depth $i$ (i.e. cliques in ${\rm TM}_i$) for $0 < i < p-1$. However this is not the case as follows from Theorem~\ref{thm:weak}: Any class of graphs which fails to be weakly $\chi_p$-bounded has already unbounded $\chi_2$ in terms of cliques at depth $1$. More formally this may be formulated as follows:

Say that a hereditary class $\mathscr C$ is {\em rank-$i$ $\chi_p$-bounded} if, for every integer $p$ there is a function $f_p$ with the property that every graph $G\in\mathscr C$ satisfies
\[
	\chi_p(G)\leq f_p(\omega({\rm TM}_{\min(i,p-1)}(G))).
\]
Then we have the following noticeable collapse.

\begin{absolutelynopagebreak}
\begin{proposition}
	A hereditary class $\mathscr C$ is rank-$i$ $\chi_p$-bounded  if and only if
\begin{itemize}
	\item $\mathscr C$ is strongly $\chi_p$-bounded and $i=0$, or
	\item $\mathscr C$ is weakly $\chi_p$-bounded and $i>0$.
\end{itemize}
\end{proposition}
\end{absolutelynopagebreak}
\begin{proof}
	If $i=0$ then the definition of a rank-$i$ $\chi_p$-bounded class reduces to the definition of a strongly $\chi_p$-bounded class.
	
	Assume $i>0$. Let $\mathscr C$ be  rank-$i$ $\chi_p$-bounded.
	Without loss of generality we can assume that all the functions $f_p$ witnessing the rank-$i$ $\chi_p$-boundedness of $\mathscr C$ are non-decreasing. Then $\mathscr C$ is rank-$(p-1)$ $\chi_p$-bounded (as ${\rm TM}_{\min(i,p-1)}(G)\subseteq {\rm TM}_{p-1}(G)$ thus $\omega({\rm TM}_{\min(i,p-1)}(G))\leq \omega({\rm TM}_{p-1}(G))$), that is $\mathscr C$ is weakly $\chi_p$-bounded. Conversely, if $\mathscr C$ is weakly $\chi_p$-bounded and, according to Theorem~\ref{thm:weak}, for every integer $s$ the class $\{G\in\mathscr C\mid K_{s,s}\not\subseteq G\}$ is strongly $\chi_p$-bounded. In other words, for every pair of integers $s$ and $t$, the class $\mathscr C_{s,t}=\{G\in\mathscr C\mid K_{s,s}\not\subseteq G\text{ and }K_t\not\subseteq G\}$ has bounded expansion. In particular there exist functions $g_p$ such that $\chi_p(G)\leq g_p(s,t)$ holds for every $G\in\mathscr C$ and $p>1$. As $s\leq \omega({\rm TM}_1(G))/2$ we deduce $\chi_p(G)\leq g_p(\omega({\rm TM}_1(G))/2,\omega({\rm TM}_1(G)))$ and thus, as $\mathscr C$ is $\chi$-bounded, the class $\mathscr C$ is rank-$1$ $\chi_p$-bounded. Hence $\mathscr C$ is rank-$i$ $\chi_p$-bounded for each $i>0$.
\myqed\end{proof}

Lemma~\ref{lem:chichis} yields that a graph with large $\chi_2$ contains, as a subgraph, the $1$-subdivision of a graph with large chromatic number. In view of Lemma~\ref{lem:chis-sub} it is perhaps natural to ask whether this property extends to all $\chi_p$.
\begin{problem}
	Is it true that for every positive integer $p$ there 
	is a function $F_p$ such that for every graph $G$ we have
\[
	\chi_p(G)\leq F_p(\chi({\rm TM}_{p-1}(G)))?
\]
	\end{problem}
\section*{Acknowledgements}
The authors would like to thank the anonymous referees for his
corrections and suggestions.

\end{document}